\documentclass[10pt]{article} 
\usepackage[left=1in,top=1in,right=1in,bottom=1in]{geometry}
\usepackage{amsfonts,amssymb,amsmath,latexsym,indentfirst,} 
 
\newtheorem{rema}{Remark}[section] 
\newtheorem{defi}{Definition}[section] 
\newtheorem{lemm}{Lemma}[section] 
\newtheorem{theo}{Theorem}[section] 
\newtheorem{coro}{Corollary}[section]

\newcommand{\R}{\ensuremath{{\mathbb{R}} }}

\newcommand{\peq}{\hspace*{0.10in}}
\newcommand{\ppeq}{\hspace*{0.05in}}
\newcommand{\para}{\hspace*{0.25in}}

\newcommand{\fim}{\rightline{$\blacksquare$}}
\newcommand{\e}{\varepsilon}

\author{Luiz Gustavo Farah \footnote{Partially supported by CNPq-Brazil.}.\\
\\
Instituto de Matem\'atica Pura e Aplicada (IMPA)\\
Estrada Dona Castorina 110, Rio de Janeiro\\
 22460-320, BRAZIL}
\title{Local solutions in Sobolev spaces with negative indices for the ``good'' Boussinesq equation } 
\date{}

\begin{document} 
\maketitle 

\begin{abstract} 
We study the local well-posedness of the initial-value problem for the nonlinear ``good'' Boussinesq equation with data in Sobolev spaces \textit{$H^s$} for negative indices of $s$.
\end{abstract} 

\section{Introduction} 
In this work we consider initial value problem (IVP) for the Boussinesq-type Equation
\begin{eqnarray}\label{NLB}
\left\{ 
\begin{array}{l}
u_{tt}-u_{xx}+u_{xxxx}+(f(u))_{xx}=0, \peq  x\in \R, t>0,\\
u(0,x)=u_0(x); \ppeq u_t(0,x)=u_1(x).
\end{array} \right. 
\end{eqnarray}

Equations of this type, but with the opposite sign in front of the fourth derivative term, were originally derived by Boussinesq \cite{BOU} in his study of nonlinear, dispersive wave propagation.
We should remark that it was the first equation proposed in the literature to describe this kind of physical phenomena. The equation (\ref{NLB}) was also used by Zakharov \cite{Z} as a model of nonlinear string and by Falk \textit{et al} \cite{FLS} in their study of shape-memory alloys.

Our principal aim here is to study the local well-posedness for the IVP associated to the ``good'' Boussinesq equation, that is, $f(u)=u^2$ in equation (\ref{NLB}), for low regularity data. Natural spaces to measure this regularity are the classical Sobolev spaces $H^s(\R)$, $s\in \R$, which are defined as the completion of the Schwarz class $\mathcal{S}(\R)$ with respect to the norm
\begin{equation*}
\|f\|_{H^{s}(\R)}=\|\langle\xi\rangle^s \widehat{f}\|_{L^{2}(\R)}
\end{equation*}
where $\langle a\rangle\equiv 1+|a|$.

Concerning the local well-posedness question, several results are obtained for the equation (\ref{NLB}).  Using Kato's abstract theory for quasilinear evolution equation, Bona and Sachs \cite{BS} showed local well-posedness for $f\in C^{\infty}$ and initial data $u_0 \in H^{s+2}(\R)$, $u_1 \in H^{s+1}(\R)$ with $s > \frac{1}{2}$ . Tsutsumi and Matahashi \cite{TM} established similar result when $f(u)=|u|^{p-1}u$, $p>1$ and $u_0 \in H^{1}(\R)$, $u_1 = \chi_{xx}$ with $\chi \in H^{1}(\R)$. These results were significantly improved by Linares \cite{FL} who proved that (\ref{NLB}) is locally well-posedness in the case $f(u)=|u|^{p-1}u$, $1<p<5$ and $u_0 \in L^{2}(\R)$, $u_1 = h_{x}$ with $h \in H^{-1}(\R)$. The main tool used in his argument was the use of Strichartz estimates satisfied by solutions of the linear problem. We should remark that all these results also hold for the ``good'' Boussinesq equation.


In this paper, we improve the latter result, proving local well-posedness with $s> -1/4$ for the ``good'' Boussinesq equation. The local well-posedness for dispersive equations with quadratic nonlinearities has been extensively studied in Sobolev spaces with negative indices. The proof of these results are based in the Fourier restriction norm approach introduced by Bourgain \cite{B} in his study of the nonlinear Schr\"odinger equation (NLS)
\begin{equation*}
iu_{t}+u_{xx}+u|u|^{p-2}=0, \textrm{ with } p\geq 3
\end{equation*}
and the Korteweg-de Vries equation (KdV)
\begin{equation}\label{KdV}
u_{t}+u_{xxx}+u_{x}u=0. 
\end{equation}
This method was further developed by Kenig, Ponce and Vega in \cite{KPV1} for the KdV equation (\ref{KdV}) and \cite{KPV2} for the quadratics for the quadratics nonlinear Schr\"odinger equations
\begin{eqnarray}\label{QNLS1}
iu_{t}+u_{xx}+F_j(u,\bar{u})=0, \,\,\, j=1,2,3
\end{eqnarray}
where $\bar{u}$ denotes the complex conjugate of $u$ and $F_1(u,\bar{u})=u^2$, $F_2(u,\bar{u})=u\bar{u}$, $F_3(u,\bar{u})={\bar{u}}^2$ in one spatial dimension and in spatially continuous and periodic case.  

The original Bourgain method makes extensive use of the Strichartz inequalities in order to derive the bilinear estimates corresponding to the nonlinearity. On the other hand, Kenig, Ponce and Vega simplified Bourgain's proof and improved the bilinear estimates using only elementary techniques, such as Cauchy-Schwartz inequality and simple calculus inequalities.

Both arguments also use some arithmetic facts involving the symbol of the  linearized equation. For example, the algebraic relation for quadratic NLS (\ref{QNLS1}) with $j=1$ is given by
\begin{equation}\label{AR}
2|\xi_1(\xi-\xi_1)|\leq |\tau-\xi^2|+|(\tau-\tau_1)-(\xi-\xi_1)^2|+|\tau_1-\xi_1^2|.
\end{equation}

Then splitting the domain of integration in the sets where each term on the right side of (\ref{AR}) is the biggest one, Kenig, Ponce and Vega made some cancellation in the symbol in order to use his calculus inequalities (see Lemma \ref{l3.1}) and a clever change of variables to established their crucial estimates. 

Here, we shall use this kind of argument, but unfortunately in the Boussinesq case we do not have good cancellations on the symbol. To overcome this difficulty we observe that the dispersion in the Boussinesq case is given by the symbol $\sqrt{{\xi}^2+{\xi}^4}$ and this is in some sense related with the Schr\"odinger symbol (see Lemma \ref{l3.3} below). Therefore, we can modify the symbols and work only with the algebraic relations for the Schr\"odinger equation already used in Kenig, Ponce and Vega \cite{KPV2} in order to derive our relevant bilinear estimates. We should remark that in the present case we have to estimate all the possible cases for the relation $\tau\pm \xi^2$ and not only the cases treated in Kenig, Ponce and Vega \cite{KPV2}.

To describe our results we define next the $X_{s,b}$ spaces related to our problem. These spaces, with $b=\frac{1}{2}$, were first defined by Fang and Grillakis \cite{FG} for the Boussinesq-type equations in the periodic case. Using these spaces and following Bourgain's argument introduced in \cite{B} they proved local well-posedness for (\ref{NLB}) with the spatial variable in the unit circle assuming $u_0 \in H^{s}$, $u_1 \in H^{-2+s}$, with $0\leq s\leq 1$ and  $|f(u)|\leq c|u|^{p}$, with $1<p<\frac{3-2s}{1-2s}$ if $0\leq s<\frac{1}{2}$ and $1<p<\infty$ if $\frac{1}{2}\leq s\leq 1$. Moreover, if $u_0 \in H^{1}$, $u_1 \in H^{-1}$ and  $f(u)= \lambda|u|^{q-1}u-|u|^{p-1}u$, with $1<q<p$ and $\lambda \in \R$ then the solution is global.

Next we give the precise definition of the $X_{s,b}$ spaces for the Boussinesq-type equation in the continuous case.

\begin{defi}\label{GAM}
For $s,b \in \R$, $X_{s,b}$ denotes the completion of the Schwartz class $\mathcal{S}(\R^2)$ with respect to the norm
\begin{equation*}
\|F\|_{X_{s,b}}=\|\langle|\tau|-\gamma(\xi)\rangle^b\langle\xi\rangle^s \widetilde{F}\|_{L^{2}_{\tau,\xi}}
\end{equation*}
where $\gamma(\xi)\equiv\sqrt{{\xi}^2+{\xi}^4}$ and $\sim$ denotes the time-space Fourier transform. 
\end{defi}

We will also need the localized $X_{s,b}$ spaces defined as follows

\begin{defi}\label{BL}
For $s,b \in \R$ and $T \geq 0$, $X_{s,b}^T$ denotes the space endowed with the norm
\begin{equation*}
\|u\|_{X_{s,b}^T}=\inf_{w\in X_{s,b}}\left\{\|w\|_{X_{s,b}}:w(t)=u(t) \textrm{ on }  [0,T]\right\}.
\end{equation*} 
\end{defi}

Now we state the main results of this paper.

\begin{theo}\label{t1.1}
Let $s >-1/4$ and $u,v\in X_{s,-a}$. Then, there exists $c>0$ such that
\begin{equation}\label{BE}
\left\|
\left(\dfrac{|\xi|^2\widetilde{uv}(\tau,\xi)}{2i\gamma(\xi)}\right)^{\sim^{-1}}
\right\|_{X_{s,-a}}\leq c\left\|u\right\|_{X_{s,b}}\left\|v\right\|_{X_{s,b}},
\end{equation}
where $\sim^{-1}$ denotes the inverse time-space Fourier transform, holds in the following cases
\begin{enumerate}
\item [($i$)] $s\geq 0$, $b>1/2$ and $1/4<a<1/2$,
\item [($ii$)]$-1/4<s<0$, $b>1/2$ and $1/4<a<1/2$ such that $|s|<a/2$.
\end{enumerate}

Moreover, the constant $c>0$ that appears in (\ref{BE}) depends only on $a,b,s$.
\end{theo}

\begin{theo}\label{t1.2}
For any $s\leq -1/4$ and any $a, b\in\R$, with $a<1/2$ the estimate (\ref{BE}) fails. 
\end{theo}

\begin{theo}\label{t1.3}
Let $s>-1/4$, then for all $\phi\in H^s(\R)$ and $\psi\in H^{s-1}(\R)$, there exist $T=T(\|\phi\|_{H^s},\|\psi\|_{H^{s-1}})$ and a unique solution $u$ of the IVP (\ref{NLB}) with $f(u)=u^2$, $u_0=\phi$ and $u_1=\psi_x$ such that
\begin{equation*}
u\in C([0,T]:H^s(\R))\cap X_{s,b}^T.
\end{equation*}

Moreover, given $T'\in (0,T)$ there exists $R=R(T')>0$ such that giving the set $ W\equiv\{(\tilde{\phi},\tilde{\psi})\in H^s(\R)\times H^{s-1}(\R):\|\tilde{\phi}-\phi\|_{H^s(\R)}^2+ \|\tilde{\psi}-\psi\|_{H^{s-1}(\R)}^2<R\}$ the map solution 
\begin{equation*}
S:W \longrightarrow C([0,T']:H^s(\R))\cap X_{s,b}^T, \peq (\tilde{\phi},\tilde{\psi})\longmapsto u(t)
\end{equation*}
is Lipschitz.

In addition, if $(\phi,\psi)\in H^{s'}(\R)\times H^{s'-1}(\R)$ with $s'>s$, then the above results hold with $s'$ instead of $s$ in the same interval $[0,T]$ with 
\begin{equation*}
T=T(\|\phi\|_{H^s},\|\psi\|_{H^{s-1}}).
\end{equation*} 
\end{theo}

Since scaling argument cannot be applied to the Boussinesq-type equations to obtain a critically notion it is not clear what is the lower index $s$ where one has local well-posedness for the ``good'' Boussinesq equation (\ref{NLB}) with with $f(u)=u^2$, $u_0=\phi$ and $u_1=\psi_x$ where $(\phi,\psi)\in H^{s}(\R)\times H^{s-1}(\R)$. Here we answer, partially, this question. In fact, our main result is a negative one; it concerns in particular a kind of ill-posedness. We prove that the flow map for the Cauchy problem (\ref{NLB}) is not smooth ($C^2$) at the origin for initial data in $H^s(\R)\times H^{s-1}(\R)$, with $s<-2$. Therefore any iterative method applied to the integral formulation of ``good'' Boussinesq equation (\ref{NLB}) always fails in this functional setting. In other words, if one can apply the contraction mapping principle to solve the integral equation corresponding to (\ref{NLB}) thus, by the implicit function Theorem, the flow-map data solution is smooth, which is a contradiction (cf. Theorem \ref{t4.2}).

Tzvetkov \cite{T} (see also Bourgain \cite{B2}) established a similar result for the KdV equation (\ref{KdV}). The same question was studied by Molinet, Saut and Tzvetkov \cite{MST1}-\cite{MST2}, for the Benjamin-Ono (BO) equation 
\begin{equation}\label{BO}
u_{t}+\mathcal{H}u_{xx}+uu_{x}=0
\end{equation}
and for the Kadomtsev-Petviashvili 1 (KP1) equation 
\begin{equation}\label{KP1}
(u_{t}+uu_{x}+u_{xxx})_x-u_{yy}=0,
\end{equation}
respectively.

Before stating the main results let us define the flow-map data solution as
\begin{equation}\label{DM}
\left. 
\begin{array}{c c c c}
S:&H^{s}(\R)\times H^{s-1}(\R)&\rightarrow&C([0,T]:H^s(\R))\\
&(\phi,\psi)&\mapsto&u(t)
\end{array} \right.
\end{equation}
where $u(t)$ is given in (\ref{INT}) below.

These are our ill-posed results
\begin{theo}\label{t4.1}
Let $s<-2$ and any $T>0$. Then there does not exist any space $X_T$ such that

\begin{equation}\label{i}
\left\|u\right\|_{C([0,T]:H^s(\R))}\leq c\left\|u\right\|_{X_T},
\end{equation}
for all $u\in X_T$
\begin{equation}\label{ii}
\left\|V_c(t)\phi+V_s(t)\psi_x\right\|_{X_T}\leq c\left(\left\|\phi\right\|_{H^s(\R)} +\left\|\psi\right\|_{H^{s-1}(\R)}\right), 
\end{equation}
for all $\phi\in H^s(\R)$, $\psi\in H^{s-1}(\R)$ and
\begin{equation}\label{iii}
\left\|\int_{0}^{t}V_s(t-t')(uv)_{xx}(t')dt'\right\|_{X_T}\leq c\left\|u\right\|_{X_T}\left\|v\right\|_{X_T},
\end{equation}
for all $u,v\in X_T$.
\end{theo}

\begin{rema}
We recall that in Sections 2 and 3 we construct a space $X_{s,b}$ such that the inequalities (\ref{i}), (\ref{ii}) and (\ref{iii}) hold for $s>-1/4$. These are the main tools to prove the local well-posedness result stated in Theorem \ref{t1.3}.
\end{rema}

\begin{theo}\label{t4.2}
Let $s<-2$. If there exists some $T>0$ such that (\ref{NLB}) with $f(u)=u^2$, $u_0=\phi$ and $u_1=\psi_x$ is locally well-posed, then the flow-map data solution $S$ defined in (\ref{DM}) is not $C^2$ at zero.
\end{theo}

In all the ill-posedness results of Tzvetkov \cite{T}, Molinet, Saut and Tzvetkov \cite{MST1}-\cite{MST2} it is, in fact, proved that for a fixed $t>0$ the flow map $S_t: \phi \mapsto u(t)$ is not $C^2$ differentiable at zero. This, of course, implies that the flow map $S$ is not smooth ($C^2$) at the origin.

Unfortunately, in our case we cannot fix $t>0$ since we don't have good cancellations on the symbol $\sqrt{{\xi}^2+{\xi}^4}$. To overcome this difficulty, we allow the variable $t$ to move. Therefore, choosing suitable characteristics functions and sending $t$ to zero we can establish Theorems \ref{t4.1}-\ref{t4.2}. We should remark that this kind of argument also appears in the ill-posed result of Bejenaru, Tao \cite{BT}. 

The plan of this paper is as follows: in Section 2, we prove some estimates for the integral equation in the $X_{s,b}$ space introduced above. Bilinear estimates and the relevant counterexamples are proved in Section 3 and 4, respectively. In Section 5, we prove the local well-posedness result stated in Theorem \ref{t1.3}. Finally, the ill-posedness question is treated in Section 6.


\section{Preliminary Results}

First, we consider the linear equation
\begin{equation}\label{LB} 
u_{tt}-u_{xx}+u_{xxxx}=0
\end{equation}

the solution for initial data $u(0)=\phi$ and $u_t(0)=\psi_x$, is given by
\begin{equation}\label{GUB} 
u(t)=V_c(t)\phi+V_s(t)\psi_x
\end{equation}

where
\begin{eqnarray*}
V_c(t)\phi&=&\left( \frac{e^{it\sqrt{{\xi}^2+{\xi}^4}}+ e^{-it\sqrt{{\xi}^2+{\xi}^4}}}{2}\hat{\phi}(\xi)\right)^{\vee}\\
V_s(t){\psi_x}&=&\left( \frac{e^{it\sqrt{{\xi}^2+{\xi}^4}}- e^{-it\sqrt{{\xi}^2+{\xi}^4}}}{2i\sqrt{{\xi}^2+{\xi}^4}}\hat{\psi_x}(\xi)\right)^{\vee}.
\end{eqnarray*}

By Duhamel's Principle the solution of (NLB) is equivalent to
\begin{equation}\label{INT} 
u(t)= V_c(t)\phi+V_s(t)\psi_x+\int_{0}^{t}V_s(t-t')(u^2)_{xx}(t')dt'.
\end{equation}

Let $\theta$ be a cutoff function satisfying $\theta \in C^{\infty}_{0}(\R)$, $0\leq \theta \leq 1$, $\theta \equiv 1$ in $[-1,1]$, supp$(\theta) \subseteq [-2,2]$ and for $0<T<1$ define $\theta_T(t)=\theta(t/T)$. In fact, to work in the $X_{s,b}$ spaces we consider another version of (\ref{INT}), that is
\begin{equation}\label{INT2} 
u(t)= \theta(t)\left(V_c(t)\phi+V_s(t)\psi_x\right)+\theta_T(t)\int_{0}^{t} V_s(t-t')(u^2)_{xx}(t')dt'.
\end{equation}

Note that the integral equation (\ref{INT2}) is defined for all $(t,x)\in \R^2$. Moreover if $u$ is a solution of (\ref{INT2}) than $\tilde{u}=u|_{[0,T]}$ will be a solution of (\ref{INT}) in $[0,T]$.\\

In the next lemma, we estimate the linear part of the integral equation (\ref{INT2}).

\begin{lemm}\label{l21}
Let $u(t)$ the solution of the linear equation 
\begin{eqnarray*}
\left\{ 
\begin{array}{l}
u_{tt}-u_{xx}+u_{xxxx}=0,\\
u(0,x)=\phi(x); \peq u_t(0,x)=(\psi(x))_x
\end{array} \right. 
\end{eqnarray*}
with $\phi \in H^s$ and $\psi \in H^{s-1}$. Then there exists $c>0$ depending only on $\theta,s,b$ such that
\begin{equation}\label{LP}
\|\theta u\|_{X_{s,b}}\leq c\left(\|\phi\|_{H^s}+\|\psi\|_{H^{s-1}}\right).
\end{equation}
\end{lemm}

\textbf{Proof. } Taking time-space Fourier transform in $\theta(t)u(x,t)$ and setting $\gamma(\xi)=\sqrt{{\xi}^2+{\xi}^4}$, we have
\begin{eqnarray*}
(\theta(t)u(x,t))^{\sim}(\tau,\xi)&=&\frac{\hat{\theta}(\tau-\gamma(\xi))}{2}
\left(\hat{\phi}(\xi)+\frac{\xi\hat{\psi}(\xi)}{\gamma(\xi)}\right)\\
&&+\frac{\hat{\theta}(\tau+\gamma(\xi))}{2}
\left(\hat{\phi}(\xi)-\frac{\xi\hat{\psi}(\xi)}{\gamma(\xi)}\right).
\end{eqnarray*}

Thus, setting $h_1(\xi)=\hat{\phi}(\xi)+\frac{\xi\hat{\psi}(\xi)}{\gamma(\xi)}$ and $h_2(\xi)=\hat{\phi}(\xi)-\frac{\xi\hat{\psi}(\xi)}{\gamma(\xi)}$, we have
\begin{eqnarray*}
\|\theta u\|_{X_{s,b}}^2\leq
\end{eqnarray*}
\begin{eqnarray*}
&\leq\int_{-\infty}^{+\infty}\langle\xi\rangle^{2s}
|h_1(\xi)|^2
\left(\int_{-\infty}^{+\infty}\langle|\tau|-\gamma(\xi)\rangle^{2b}
\left|\frac{\hat{\theta}(\tau-\gamma(\xi))+\hat{\theta}(\tau+\gamma(\xi))}{2}\right|^2
d\tau\right)d\xi\\
&+\int_{-\infty}^{+\infty}\langle\xi\rangle^{2s}
|h_2(\xi)|^2
\left(\int_{-\infty}^{+\infty}\langle|\tau|-\gamma(\xi)\rangle^{2b}
\left|\frac{\hat{\theta}(\tau-\gamma(\xi))+\hat{\theta}(\tau+\gamma(\xi))}{2}\right|^2
d\tau\right)d\xi.
\end{eqnarray*}

Since $||\tau|-\gamma(\xi)|\leq \min\left\{|\tau-\gamma(\xi)|, |\tau+\gamma(\xi)|\right\}$ and $\hat{\theta}$ is rapidly decreasing, we can bound the terms inside the parentheses, and the claim follows. 
\fim

Next we estimate the integral part of (\ref{INT2}).

\begin{lemm}\label{L22}
Let $-\frac{1}{2}<b'\leq 0\leq b \leq b'+1$ and $0<T \leq 1$ then
\begin{enumerate}
\item [$(i)$] $\left\|\theta_T(t)\int_{0}^{t}g(t')dt'\right\|_{H^b_t} \leq T^{1-(b-b')}\|g\|_{H^{b'}_{t}}$;

\item [$(ii)$] $\left\|\theta_T(t)\int_{0}^{t}V_s(t-t')f(u)(t')dt'\right\|_{X_{s,b}} \leq T^{1-(b-b')} \left\|
\left(\dfrac{\widetilde{f(u)}(\tau,\xi)}{2i\gamma(\xi)}\right)^{\sim^{-1}}
\right\|_{X_{s,b'}}$.
\end{enumerate}
\end{lemm}

\textbf{Proof. } 
\begin{enumerate}
\item [$(i)$] See \cite{G} inequality $(3.11)$.
\item [$(ii)$] A simple calculation shows that
\begin{eqnarray*}
\left(\theta_T(t)\int_{0}^{t}V_s(t-t')f(u)(t') dt'\right)^{\wedge_{(x)}}(t,\xi)=
\end{eqnarray*}
\begin{eqnarray*}
&=&e^{it\gamma(\xi)} \left(\theta_T(t)\int_{0}^{t}h_1(t',\xi)dt'\right)
-e^{-it\gamma(\xi)} \left(\theta_T(t)\int_{0}^{t}h_2(t',\xi)dt'\right)\\
&\equiv& e^{it\gamma(\xi)}w^{\wedge_{(x)}}_1(t,\xi)-e^{-it\gamma(\xi)}w^{\wedge_{(x)}}_2(t,\xi),\\
\end{eqnarray*}
where $h_1(t',\xi)=\dfrac{e^{-it'\gamma(\xi)}f^{\wedge_{(x)}}(t',\xi)}{2i\gamma(\xi)}$ and $h_2(t',\xi)=\dfrac{e^{it'\gamma(\xi)}f^{\wedge_{(x)}}(t',\xi)}{2i\gamma(\xi)}$.
\end{enumerate}

Therefore
\begin{eqnarray*}
\left(\theta_T(t)\int_{0}^{t}V_s(t-t')f(u)(t') dt'\right)^{\sim}(\tau,\xi)=
\end{eqnarray*}
\begin{eqnarray*}
\widetilde{w_1}(\tau-\gamma(\xi),\xi)- \widetilde{w_2}(\tau+\gamma(\xi),\xi).
\end{eqnarray*}

Now using the definition of $X_{s,b}$ we have
\begin{eqnarray*}
\left\|\theta_T(t)\int_{0}^{t}V_s(t-t')f(u)(t') dt'\right\|_{X_{s,b}}^2\leq
\end{eqnarray*}
\begin{eqnarray*}
&\leq&\int_{-\infty}^{+\infty}\int_{-\infty}^{+\infty} \langle|\tau+\gamma(\xi)|-\gamma(\xi)\rangle^{2b}\langle\xi\rangle^{2s} |\widetilde{w_1}(\tau,\xi)|^2d\tau d\xi\\
&&+\int_{-\infty}^{+\infty}\int_{-\infty}^{+\infty} \langle|\tau-\gamma(\xi)|-\gamma(\xi)\rangle^{2b}\langle\xi\rangle^{2s} |\widetilde{w_2}(\tau,\xi)|^2d\tau d\xi\\
&\equiv& M.
\end{eqnarray*}

Since $\gamma(\xi)\geq 0$ for all $\xi \in \R$, we have
\begin{eqnarray*}
\max\{||\tau+\gamma(\xi)|-\gamma(\xi)|,||\tau-\gamma(\xi)|-\gamma(\xi)|\}\leq |\tau|.
\end{eqnarray*}

Thus applying item $(i)$ we obtain
\begin{eqnarray*}
M&\leq&c\sum_{j=1}^{2}\int_{-\infty}^{+\infty}\langle\xi\rangle^{2s} \|w_j^{\wedge_{(x)}}\|_{H^b_t}^2\\
&\leq&cT^{1-(b-b')}\left(\int\int_{\R^2}\langle\tau-\gamma(\xi)\rangle^{2b'}\langle\xi\rangle^{2s} \left|\dfrac{\widetilde{f(u)}(\tau,\xi)}{2i\gamma(\xi)}\right|^2d\tau d\xi\right.\\
&&+\left.\int\int_{\R^2}\langle\tau+\gamma(\xi)\rangle^{2b'}\langle\xi\rangle^{2s} \left|\dfrac{\widetilde{f(u)}(\tau,\xi)}{2i\gamma(\xi)}\right|^2d\tau d\xi\right).
\end{eqnarray*}

Since $||\tau|-\gamma(\xi)|\leq \min\left\{|\tau-\gamma(\xi)|, |\tau+\gamma(\xi)|\right\}$ and $b'\leq 0$ we obtain the desired inequality.\\
\fim

The next lemma says that, for $b> 1/2$, $X_{s,b}$ is embedding in $C(\R:H^s)$.

\begin{lemm}\label{l11}
Let $b>\frac{1}{2}$. There exists $c>0$, depending only on $b$, such that
\begin{equation*}
\|u\|_{C(\R:H^s)}\leq c\|u\|_{X_{s,b}}.
\end{equation*}
\end{lemm}

\textbf{Proof. } First we prove that $X_{s,b}\subseteq L^{\infty}(\R:H^s)$. Let $u = u_1+u_2$, where $\tilde{u}_1\equiv\tilde{u}\chi_{\{\tau \leq 0\}}$, $\tilde{u}_2\equiv\tilde{u}\chi_{\{\tau > 0\}}$ and $\chi_A$ denotes the characteristic function of the set $A$. Then for all $t\in \R$
\begin{eqnarray*}
\|u_1(t)\|_{H^s}&=&\left\|\left( e^{it\gamma(\xi)}(u_1)^{\wedge_{(x)}}\right)^{\vee_{(x)}}(t,x)\right\|_{H^s}\\
&=&\left\|\int_{-\infty}^{+\infty} \left(\left(e^{it\gamma(\xi)}(u_1)^{\wedge_{(x)}}\right)^{\vee_{(x)}} \right)^{\wedge_{(t)}} (\tau,x)e^{it\tau}d\tau\right\|_{H^s}\\
&\leq&\int_{-\infty}^{+\infty}\left\| \left(\left(e^{it\gamma(\xi)}(u_1)^{\wedge_{(x)}}\right)^{\vee_{(x)}}\right) ^{\wedge_{(t)}}(\tau,x) \right\|_{H^s}d\tau.\\
\end{eqnarray*}

Using the Cauchy-Schwarz inequality we obtain
\begin{eqnarray*}
\|u_1(t)\|_{H^s}
\leq&\left(\int_{-\infty}^{+\infty}\langle\tau\rangle^{-2b}\right)^{1/2} \left(\int_{-\infty}^{+\infty}\int_{-\infty}^{0} \langle\tau+\gamma(\xi)\rangle^{2b}\langle\xi\rangle^{2s} |\tilde{u}(\tau,\xi)|^2d\tau d\xi\right)^{1/2}.
\end{eqnarray*}

On the other hand, similar arguments imply that
\begin{eqnarray*}
\|u_2(t)\|_{H^s}
\leq&\left(\int_{-\infty}^{+\infty}\langle\tau\rangle^{-2b}\right)^{1/2} \left(\int_{-\infty}^{+\infty}\int_{0}^{+\infty} \langle\tau-\gamma(\xi)\rangle^{2b}\langle\xi\rangle^{2s} |\tilde{u}(\tau,\xi)|^2d\tau d\xi\right)^{1/2}.
\end{eqnarray*}
 
Now, by the fact that $b>1/2$, $|\tau+\gamma(\xi)|=||\tau|-\gamma(\xi)|$ for $\tau \leq 0$ and $|\tau-\gamma(\xi)|=||\tau|-\gamma(\xi)|$ for $\tau \geq 0$ we have
\begin{equation*}
\|u\|_{L^{\infty}(\R:H^s)}\leq c\|u\|_{X_{s,b}}.
\end{equation*}

It remains to show continuity. Let $t,t'\in \R$ then
\begin{eqnarray*}
\|u_1(t)-u_1(t')\|_{H^s}=
\end{eqnarray*}
\begin{eqnarray}\label{TCO}
\left\|\int_{-\infty}^{+\infty}\left(\left(e^{it\gamma(\xi)}(u_1)^{\wedge_{(x)}}\right)^{\vee_{(x)}} \right)^{\wedge_{(t)}} (\tau,x)(e^{it\tau}-e^{it'\tau})d\tau\right\|_{H^s}
\end{eqnarray}

Letting $t'\rightarrow t$, two applications of the Dominated Convergence Theorem give that the right hand side of (\ref{TCO}) goes to zero. Therefore, $u_1\in C(\R:H^s)$. Of course, the same argument applies to $u_2$, which concludes the proof.
\fim

To finish this section, we remark that for any positive numbers $a$ and $b$, the notation $a \lesssim b$ means that there exists a positive
constant $\theta$ such that $a \leq \theta b$. We also denote $a \sim b$ when, $a \lesssim b$ and $b \lesssim a$.

\section{Bilinear estimates}

Before proceed to the proof of Theorem \ref{t1.1}, we state some elementary calculus inequalities that will be useful later.
\begin{lemm}\label{l3.1}
For $p, q>0$ and $r=\min\{ p, q, p+q-1\}$ with $p+q>1$, there exists $c>0$ such that
\begin{equation}\label{CI1}
\int_{-\infty}^{+\infty}\dfrac{dx} {\langle x-\alpha\rangle^{p}\langle x-\beta\rangle^{q}}\leq\dfrac{c} {\langle \alpha-\beta\rangle^{r}}.
\end{equation}
Moreover, for $a_0, a_1, a_2\in \R$ and $q>1/2$
\begin{equation}\label{CI2}
\int_{-\infty}^{+\infty}\dfrac{dx} {\langle a_0+a_1x+a_2x^2\rangle^{q}}\leq c. 
\end{equation}
\end{lemm}
\textbf{Proof. } See Lemma 4.2 in \cite{GTV} and Lemma 2.5 in \cite{BOP}.\\
\fim

\begin{lemm}\label{l3.3}
There exists $c>0$ such that
\begin{equation}\label{LN}
\dfrac{1}{c}\leq\sup_{x,y\geq 0}\dfrac{1+|x-y|}{1+|x-\sqrt{y^2+y}|}\leq c.
\end{equation}
\end{lemm}

\textbf{Proof. } Since $y\leq\sqrt{y^2+y}\leq y+1/2$ for all $y\geq 0$ a simple computation shows the desired inequalities.\\
\fim
\begin{rema}\label{R1}
In view of the previous lemma we have an equivalent way to compute the $X_{s,b}$-norm, that is 
\begin{equation*}
\|u\|_{X_{s,b}}\sim\|\langle|\tau|-\xi^2\rangle^b\langle\xi\rangle^{s} \widetilde{u}(\tau,\xi)\|_{L^2_{\tau,\xi}}. 
\end{equation*}
This equivalence will be important in the proof of Theorem \ref{t1.1}. As we commented in the introduction the Boussinesq symbol $\sqrt{{\xi}^2+{\xi}^4}$ does not have good cancelations to make use of Lemma \ref{l3.1}. Therefore, we  modify the symbols as above and work only with the algebraic relations for the Schr\"odinger equation already used in Kenig, Ponce and Vega \cite{KPV2} in order to derive the bilinear estimates. 
\end{rema}

Now we are in position to prove the bilinear estimate (\ref{BE}).\\

\textbf{Proof of Theorem \ref{t1.1}.} Let $u, v\in X_{s,b}$ and define $f(\tau,\xi)\equiv \langle|\tau|-\xi^2\rangle^b\langle\xi\rangle^{s} \widetilde{u}(\tau,\xi)$, $g(\tau,\xi)\equiv \langle|\tau|-\xi^2\rangle^b\langle\xi\rangle^{s} \widetilde{v}(\tau,\xi)$. Using Remark \ref{R1} and a duality argument the desired inequality is equivalent to
\begin{equation}\label{DUA}
\left|W(f,g,\phi)\right|\leq c\|f\|_{L^2_{\xi,\tau}}\|g\|_{L^2_{\xi,\tau}}\|\phi\|_{L^2_{\xi,\tau}}
\end{equation}
where
\begin{eqnarray*}
W(f,g,\phi)&=&  \int_{\R^4} \dfrac{\langle\xi\rangle^{s}}{\langle\xi_1\rangle^{s} \langle\xi-\xi_1\rangle^{s} }\\
&&\times \dfrac{g(\tau_1,\xi_1)f(\tau-\tau_1,\xi-\xi_1) \bar{\phi}(\tau,\xi)}{\langle|\tau|-\xi^2\rangle^{a} \langle|\tau_1|-\xi_1^2\rangle^{b} \langle|\tau-\tau_1|-(\xi-\xi_1)^2\rangle^{b} }d\xi d\tau d\xi_1 d\tau_1
\end{eqnarray*}

Therefore to perform the desired estimate we need to analyze all the possible cases for the sign of $\tau$, $\tau_1$ and $\tau-\tau_1$. To do this we split $\R^4$ into the regions 
\begin{eqnarray*}
\Gamma_1&=&\{(\xi, \tau, \xi_1, \tau_1)\in \R^4: \tau_1< 0, \tau-\tau_1< 0\},\\ 
\Gamma_2&=&\{(\xi, \tau, \xi_1, \tau_1)\in \R^4: \tau_1\geq 0, \tau-\tau_1< 0, \tau\geq 0\},\\ 
\Gamma_3&=&\{(\xi, \tau, \xi_1, \tau_1)\in \R^4: \tau_1\geq 0, \tau-\tau_1< 0, \tau< 0\},\\ 
\Gamma_4&=&\{(\xi, \tau, \xi_1, \tau_1)\in \R^4: \tau_1< 0, \tau-\tau_1\geq 0, \tau\geq 0\},\\
\Gamma_5&=&\{(\xi, \tau, \xi_1, \tau_1)\in \R^4: \tau_1< 0, \tau-\tau_1\geq 0, \tau< 0\},\\
\Gamma_6&=&\{(\xi, \tau, \xi_1, \tau_1)\in \R^4: \tau_1\geq 0, \tau-\tau_1\geq 0\}.
\end{eqnarray*}

Thus, it is sufficient to prove inequality (\ref{DUA}) with $Z(f,g,\phi)$ instead of $W(f,g,\phi)$, where
 \begin{equation*}
Z(f,g,\phi)=  \int_{\R^4} \dfrac{\langle\xi\rangle^{s}}{\langle\xi_1\rangle^{s} \langle\xi_2\rangle^{s}} \dfrac{g(\tau_1,\xi_1)f(\tau_2,\xi_2) \bar{\phi}(\tau,\xi)}{\langle\sigma\rangle^{a} \langle\sigma_1\rangle^{b} \langle\sigma_2\rangle^{b}}d\xi d\tau d\xi_1 d\tau_1
\end{equation*}
with $\xi_2=\xi-\xi_1$, $\tau_2=\tau-\tau_1$ and $\sigma, \sigma_1, \sigma_2$ belonging to one of the following cases
\begin{enumerate}
\item [$(I)$] $\sigma=\tau+\xi^2,\peq \sigma_1=\tau_1+\xi_1^2,\peq \sigma_2=\tau_2+\xi_2^2$,
\item [$(II)$]$\sigma=\tau-\xi^2,\peq \sigma_1=\tau_1-\xi_1^2,\peq \sigma_2=\tau_2+\xi_2^2$,
\item [$(III)$] $\sigma=\tau+\xi^2,\peq \sigma_1=\tau_1-\xi_1^2,\peq \sigma_2=\tau_2+\xi_2^2$,
\item [$(IV)$] $\sigma=\tau-\xi^2,\peq \sigma_1=\tau_1+\xi_1^2,\peq \sigma_2=\tau_2-\xi_2^2$,
\item [$(V)$] $\sigma=\tau+\xi^2,\peq \sigma_1=\tau_1+\xi_1^2,\peq \sigma_2=\tau_2-\xi_2^2$,
\item [$(VI)$] $\sigma=\tau-\xi^2,\peq \sigma_1=\tau_1-\xi_1^2,\peq \sigma_2=\tau_2-\xi_2^2$.
\end{enumerate}

\begin{rema}
Note that the cases $\sigma=\tau+\xi^2,\peq \sigma_1=\tau_1-\xi_1^2,\peq \sigma_2=\tau_2-\xi_2^2$ and $\sigma=\tau-\xi^2,\peq \sigma_1=\tau_1+\xi_1^2,\peq \sigma_2=\tau_2+\xi_2^2$ cannot occur, since $\tau_1< 0, \tau-\tau_1< 0$ implies $\tau<0$ and $\tau_1\geq 0, \tau-\tau_1\geq 0$ implies $\tau\geq 0$
\end{rema}
Applying the change of variables $(\xi, \tau, \xi_1, \tau_1)\mapsto -(\xi, \tau, \xi_1, \tau_1)$ and observing that the $L^2$-norm is preserved under the reflection operation, the cases $(IV)$, $(V)$, $(VI)$ can be easily reduced, respectively, to $(III)$, $(II)$, $(I)$. Moreover,  making the change of variables $\tau_2=\tau-\tau_1$, $\xi_2=\xi-\xi_1$ and then $(\xi, \tau, \xi_2, \tau_2)\mapsto -(\xi, \tau, \xi_2, \tau_2)$ the case $(II)$ can be reduced $(III)$. Therefore we need only establish cases $(I)$ and $(III)$. We should remark that these are exactly the two estimates that appear in \cite{KPV2}, but since it is important to have the inequality (\ref{DUA}) with $a<1/2<b$ such that $a+b<1$ to make the contraction arguments work (see the proof of Theorem \ref{t1.3}) we reprove these inequalities here.

We first treat the inequality (\ref{DUA}) with $Z(f,g,\phi)$ in the case ($I$). We will make use of the following algebraic relation
\begin{equation}\label{AR1}
 -(\tau+\xi^2)+(\tau_1+\xi_1^2)+((\tau-\tau_1)+(\xi-\xi_1)^2)=2\xi_1(\xi_1-\xi).
\end{equation}

By symmetry we can restrict ourselves to the set
\begin{equation*}
 A=\{(\xi, \tau, \xi_1, \tau_1)\in \R^4: |(\tau-\tau_1)+(\xi-\xi_1)^2|\leq|\tau_1+\xi_1^2|\}.
\end{equation*}

We divide $A$ into three pieces
\begin{eqnarray*}
A_1&=&\{(\xi, \tau, \xi_1, \tau_1)\in A: |\xi_1|\leq 10\}\\
A_2&=&\{(\xi, \tau, \xi_1, \tau_1)\in A: |\xi_1|\geq 10 \textrm{ and } |2\xi_1-\xi|\geq |\xi_1|/2 \}\\
A_3&=&\{(\xi, \tau, \xi_1, \tau_1)\in A: |\xi_1|\geq 10 \textrm{ and } |\xi_1-\xi|\geq |\xi_1|/2 \}.
\end{eqnarray*} 

We have $A=A_1\cup A_2\cup A_3$. Indeed
\begin{equation*}
|\xi_1|>|2\xi_1-\xi|+|\xi_1-\xi|\geq |(2\xi_1-\xi) -(\xi_1-\xi)|=|\xi_1|.
\end{equation*}

Next we split $A_3$ into two parts
\begin{eqnarray*}
A_{3,1}&=&\{(\xi, \tau, \xi_1, \tau_1)\in A_3: |\tau_1+\xi_1^2|\leq|\tau+\xi^2|\}\\
A_{3,2}&=&\{(\xi, \tau, \xi_1, \tau_1)\in A_3: |\tau+\xi^2|\leq|\tau_1+\xi_1^2| \}.
\end{eqnarray*} 

We can now define the sets $R_i$, $i=1,2$, as follows
\begin{equation*}
R_1=A_1\cup A_2\cup A_{3,1} \textrm{  and  } R_2=A_{3,2}.
\end{equation*}

In what follows $\chi_R$ denotes the characteristic function of the set $R$. Using the Cauchy-Schwarz and H\"older inequalities it is easy to see that
\begin{eqnarray*}
|Z|^2&\leq& \|f\|_{L^2_{\xi,\tau}}^2\|g\|_{L^2_{\xi,\tau}}^2\|\phi\|_{L^2_{\xi,\tau}}^2\\
&&\times\left\|\dfrac{\langle\xi\rangle^{2s}}{\langle\sigma\rangle^{2a}} \int\!\!\!\!\int\dfrac{\chi_{R_1}d\xi_1 d\tau_1}{\langle\xi_1\rangle^{2s} \langle\xi_2\rangle^{2s} \langle\sigma_1\rangle^{2b} \langle\sigma_2\rangle^{2b}}\right\|_{L^{\infty}_{\xi,\tau}}\\
&&+\|f\|_{L^2_{\xi,\tau}}^2\|g\|_{L^2_{\xi,\tau}}^2\|\phi\|_{L^2_{\xi,\tau}}^2\\
&&\times\left\|\dfrac{1}{\langle\xi_1\rangle^{2s}\langle\sigma_1\rangle^{2b}} \int\!\!\!\!\int\dfrac{\chi_{R_2}\langle\xi\rangle^{2s}d\xi d\tau}{\langle\xi_2\rangle^{2s}\langle\sigma\rangle^{2a}   \langle\sigma_2\rangle^{2b}}\right\|_{L^{\infty}_{\xi_1,\tau_1}}.\\
\end{eqnarray*}

Noting that $\langle\xi\rangle^{2s}\leq\langle\xi_1\rangle^{2|s|}\langle\xi_2\rangle^{2s}$, for $s\geq 0$, and $\langle\xi_2\rangle^{-2s}\leq\langle\xi_1\rangle^{2|s|}\langle\xi\rangle^{-2s}$, for $s< 0$ we have
\begin{equation}\label{Xi}
\dfrac{\langle\xi\rangle^{2s}}{\langle\xi_1\rangle^{2s} \langle\xi_2\rangle^{2s}}\leq \langle\xi_1\rangle^{\gamma(s)}
\end{equation}
where 
\begin{eqnarray*}
\gamma(s)=
\left\{ 
\begin{array}{l c}
0, &\textrm{ if } s\geq 0\\
4|s|, &\textrm{ if } s\leq 0
\end{array} \right. .
\end{eqnarray*}

Therefore in view of Lemma \ref{l3.1}-(\ref{CI1}) it suffices to get bounds for
\begin{eqnarray*}
J_1(\xi,\tau)&\equiv&\dfrac{1}{\langle\sigma\rangle^{2a}} \int\dfrac{\langle\xi_1\rangle^{\gamma(s)}d\xi_1} {\langle\tau+\xi^2+2\xi_1^2-2\xi\xi_1\rangle^{2b}}\peq \textrm{on} \peq R_1\\
J_2(\xi_1,\tau_1)&\equiv&\dfrac{\langle\xi_1\rangle^{\gamma(s)}} {\langle\sigma_1\rangle^{2b}} \int\dfrac{d\xi} {\langle\tau_1-\xi_1^2+2\xi\xi_1\rangle^{2a}}\peq \textrm{on} \peq R_2.
\end{eqnarray*}

In region $A_1$ we have $\langle\xi_1\rangle^{\gamma(s)}\lesssim 1$. Therefore for $a>0$ and $b>1/2$ we obtain
\begin{equation*}
J_1(\xi,\tau) \lesssim \int_{|\xi_1|\leq 10}d\xi_1 \lesssim 1.
\end{equation*}

In region $A_2$, by the change of variables $\eta=\tau+\xi^2+2\xi_1^2-2\xi\xi_1$ and the condition $|2\xi_1-\xi|\geq |\xi_1|/2$ we have
\begin{eqnarray*}
J_1(\xi,\tau) &\lesssim& \dfrac{1}{\langle\sigma\rangle^{2a}}\int \dfrac{\langle\xi_1\rangle^{\gamma(s)}}{|2\xi_1-\xi|\langle\eta\rangle^{2b}}d\eta \\
&\lesssim& \dfrac{1}{\langle\sigma\rangle^{2a}}\int \dfrac{\langle\xi_1\rangle^{\gamma(s)-1}}{\langle\eta\rangle^{2b}}d\eta
\lesssim 1
\end{eqnarray*}
for $a>0$, $b>1/2$ and $s> -1/4$ which implies ${\gamma(s)}\leq 1$.\\

Now, by definition of region $A_{3,1}$ and the algebraic relation (\ref{AR1}) we have
\begin{equation*}
\langle\xi_1\rangle^2\lesssim|\xi_1|^2\lesssim|\xi_1(\xi_1-\xi)|\lesssim \langle\sigma\rangle.
\end{equation*} 

Therefore by Lemma \ref{l3.1}-(\ref{CI2})
\begin{eqnarray*}
J_1(\xi,\tau) &\lesssim& \int \dfrac{\langle\xi_1\rangle^{\gamma(s)-4a}}{\langle\tau+\xi^2+2\xi_1^2-2\xi\xi_1\rangle^{2b}}d\xi_1 \\
&\lesssim& \int \dfrac{1}{\langle\tau+\xi^2+2\xi_1^2-2\xi\xi_1\rangle^{2b}}d\xi_1
\lesssim 1
\end{eqnarray*}
for $a>1/4$, $b>1/2$ and $s> -1/4$ which implies $\gamma(s)<4a$.\\

Next we estimate $J_2(\xi_1,\tau_1)$. Making the change of variables, $\eta=\tau-\xi_1^2+2\xi\xi_1$, using the restriction in the region $A_{3,2}$, we have
\begin{equation*}
|\eta|\lesssim|(\tau-\tau_1)+(\xi-\xi_1)^2|+|\tau+\xi^2|\lesssim \langle\sigma_1\rangle.
\end{equation*} 

Moreover, in $A_{3,2}$
\begin{equation*}
|\xi_1|^2\lesssim|\xi_1(\xi_1-\xi)|\lesssim \langle\sigma_1\rangle.
\end{equation*} 

Therefore, since $|\xi_1|\geq 10$ we have
\begin{eqnarray*}
J_2(\xi_1,\tau_1) &\lesssim& \dfrac{|\xi_1|^{\gamma(s)}}{\langle\sigma_1\rangle^{2b}} \int_{|\eta| \lesssim \langle\sigma_1\rangle} \dfrac{d\eta} {|\xi_1|\langle\eta\rangle^{2a}}\\
&\lesssim&\dfrac{|\xi_1|^{\gamma(s)-1}}{\langle\sigma_1\rangle^{2b+2a-1}} 
\lesssim 1
\end{eqnarray*}
for $0<a<1/2$, $b>1/2$ and $s> -1/4$.\\

Now we turn to the proof of case ($III$). In the following estimates we will make use of the algebraic relation
\begin{equation}\label{AR2}
 -(\tau+\xi^2)+(\tau_1-\xi_1^2)+((\tau-\tau_1)+(\xi-\xi_1)^2)=-2\xi_1\xi.
\end{equation}

First we split $\R^4$ into four sets
\begin{eqnarray*}
B_1&=&\{(\xi, \tau, \xi_1, \tau_1)\in \R^4: |\xi_1|\leq 10\}\\
B_2&=&\{(\xi, \tau, \xi_1, \tau_1)\in \R^4: |\xi_1|\geq 10 \textrm{ and } |\xi|\leq 1 \}\\
B_3&=&\{(\xi, \tau, \xi_1, \tau_1)\in \R^4: |\xi_1|\geq 10, |\xi|\geq 1 \textrm{ and } |\xi|\geq |\xi_1|/2 \}\\
B_4&=&\{(\xi, \tau, \xi_1, \tau_1)\in \R^4: |\xi_1|\geq 10, |\xi|\geq 1 \textrm{ and } |\xi|\leq |\xi_1|/2 \}.
\end{eqnarray*} 

Next we separate $B_4$ into three parts
\begin{small}
\begin{eqnarray*}
B_{4,1}&=&\{(\xi, \tau, \xi_1, \tau_1)\in B_4: |\tau_1-\xi_1^2|,|(\tau-\tau_1)+(\xi-\xi_1)^2|\leq|\tau+\xi^2|\}\\
B_{4,2}&=&\{(\xi, \tau, \xi_1, \tau_1)\in B_4: |\tau+\xi^2|,|(\tau-\tau_1)+(\xi-\xi_1)^2|\leq|\tau_1-\xi_1^2|\}\\
B_{4,3}&=&\{(\xi, \tau, \xi_1, \tau_1)\in B_4: |\tau_1-\xi_1^2|,|\tau+\xi^2|\leq|(\tau-\tau_1)+(\xi-\xi_1)^2|\}.
\end{eqnarray*} 
\end{small}

We can now define the sets $R_i$, $i=1,2,3$, as follows
\begin{equation*}
S_1=B_1\cup B_3\cup B_{4,1},\peq S_2=B_2\cup B_{4,2} \peq\textrm{  and  } \peq S_3=B_{4,3}.
\end{equation*}

Using the Cauchy-Schwarz and H\"older inequalities and duality it is easy to see that
\begin{eqnarray*}
|Z|^2&\leq& \|f\|_{L^2_{\xi,\tau}}^2\|g\|_{L^2_{\xi,\tau}}^2\|\phi\|_{L^2_{\xi,\tau}}^2\\
&&\times\left\|\dfrac{\langle\xi\rangle^{2s}}{\langle\sigma\rangle^{2a}} \int\!\!\!\!\int\dfrac{\chi_{S_1}d\xi_1 d\tau_1}{\langle\xi_1\rangle^{2s} \langle\xi_2\rangle^{2s} \langle\sigma_1\rangle^{2b} \langle\sigma_2\rangle^{2b}}\right\|_{L^{\infty}_{\xi,\tau}}\\
&&+\|f\|_{L^2_{\xi,\tau}}^2\|g\|_{L^2_{\xi,\tau}}^2\|\phi\|_{L^2_{\xi,\tau}}^2\\
&&\times\left\|\dfrac{1}{\langle\xi_1\rangle^{2s}\langle\sigma_1\rangle^{2b}} \int\!\!\!\!\int\dfrac{\chi_{S_2}\langle\xi\rangle^{2s}d\xi d\tau}{\langle\xi_2\rangle^{2s}\langle\sigma\rangle^{2a}   \langle\sigma_2\rangle^{2b}}\right\|_{L^{\infty}_{\xi_1,\tau_1}}\\
&&+\|f\|_{L^2_{\xi,\tau}}^2\|g\|_{L^2_{\xi,\tau}}^2\|\phi\|_{L^2_{\xi,\tau}}^2\\
&&\times\left\|\dfrac{1}{\langle\xi_2\rangle^{2s}\langle\sigma_2\rangle^{2b}} \int\!\!\!\!\int\dfrac{\chi_{\widetilde{S}_3}\langle\xi_1+\xi_2\rangle^{2s}d\xi_1 d\tau_1}{\langle\xi_1\rangle^{2s}\langle\sigma_1\rangle^{2a}   \langle\sigma\rangle^{2b}}\right\|_{L^{\infty}_{\xi_2,\tau_2}}.\\
\end{eqnarray*}

where $\sigma$, $\sigma_1$, $\sigma_2$ were given in the condition ($III$) and
\begin{eqnarray*}
\widetilde{S}_3\subseteq 
\left\{ 
\begin{array}{r}
(\xi_2, \tau_2, \xi_1, \tau_1)\in \R^4: |\xi_1|\geq 10, |\xi_1+\xi_2|\geq 1,  |\xi_1+\xi_2|\leq |\xi_1|/2 \textrm{ and }\\
|\tau_1-\xi_1^2|,|(\tau_1+\tau_2)+(\xi_1+\xi_2)^2|\leq|\tau_2+\xi_2^2|
\end{array} \right\}.
\end{eqnarray*} 

Therefore from Lemma \ref{l3.1}-(\ref{CI1}) and (\ref{Xi}) it suffices to get bounds for
\begin{eqnarray*}
K_1(\xi,\tau)&\equiv&\dfrac{1}{\langle\sigma\rangle^{2a}} \int\dfrac{\langle\xi_1\rangle^{\gamma(s)}d\xi_1} {\langle\tau+\xi^2-2\xi\xi_1\rangle^{2b}}\peq \textrm{on} \peq S_1\\
K_2(\xi_1,\tau_1)&\equiv&\dfrac{\langle\xi_1\rangle^{\gamma(s)}} {\langle\sigma_1\rangle^{2b}} \int\dfrac{d\xi} {\langle\tau_1-\xi_1^2+2\xi\xi_1\rangle^{2a}}\peq \textrm{on} \peq S_2\\
K_3(\xi_1,\tau_1)&\equiv&\dfrac{1}{\langle\sigma_2\rangle^{2b}} \int\dfrac{\langle\xi_1\rangle^{\gamma(s)}d\xi_1} {\langle\tau_2+\xi_2^2+2\xi_1^2+2\xi_1\xi_2\rangle^{2a}}\peq \textrm{on} \peq \widetilde{S}_3.
\end{eqnarray*}

In region $B_1$ we have $\langle\xi_1\rangle^{\gamma(s)}\lesssim 1$. Therefore for $a>0$ and $b>1/2$ we obtain
\begin{equation*}
K_1(\xi,\tau) \lesssim \int_{|\xi_1|\leq 10}d\xi_1 \lesssim 1.
\end{equation*}

In region $B_3$, the change of variables $\eta=\tau+\xi^2-2\xi\xi_1$ and the condition $|\xi|\geq |\xi_1|/2$ give
\begin{eqnarray*}
K_1(\xi,\tau) &\lesssim& \dfrac{1}{\langle\sigma\rangle^{2a}}\int \dfrac{\langle\xi_1\rangle^{\gamma(s)}}{|\xi|\langle\eta\rangle^{2b}}d\eta \\
&\lesssim& \dfrac{\langle\xi_1\rangle^{\gamma(s)-1}}{\langle\sigma\rangle^{2a}}\int \dfrac{1}{\langle\eta\rangle^{2b}}d\eta
\lesssim 1
\end{eqnarray*}
for $a>0$, $b>1/2$ and $s> -1/4$ which implies ${\gamma(s)}\leq 1$.\\

Now, by definition of region $B_{4,1}$ and the algebraic relation (\ref{AR2}) we have
\begin{equation*}
\langle\xi_1\rangle\lesssim|\xi_1|\lesssim|\xi_1\xi|\lesssim \langle\sigma\rangle.
\end{equation*} 

Therefore the change of variables $\eta=\tau+\xi^2-2\xi\xi_1$ and the condition $|\xi|\geq 1$ yield
\begin{eqnarray*}
K_1(\xi,\tau) &\lesssim& \dfrac{1}{\langle\sigma\rangle^{2a}}\int \dfrac{\langle\xi_1\rangle^{\gamma(s)}}{|\xi|\langle\eta\rangle^{2b}}d\eta \\
&\lesssim& \dfrac{\langle\xi_1\rangle^{\gamma(s)-2a}}{|\xi|}\int \dfrac{1}{\langle\eta\rangle^{2b}}d\eta
\lesssim 1
\end{eqnarray*}
for $s>-1/4$, $b>1/2$ and $a\in \R$ such that $2|s|<a<1/2$, if $s<0$ or $0<a<1/2$, if $s\geq 0$.\\

Next we estimate $K_2(\xi_1,\tau_1)$. Making the change of variables, $\eta=\tau_1-\xi_1^2+2\xi\xi_1$, using the restriction in the region $B_2$, we have
\begin{equation*}
|\eta|\lesssim|\tau_1-\xi_1^2|+|\xi\xi_1|\lesssim |\sigma_1|+|\xi_1|.
\end{equation*} 

Therefore, 
\begin{eqnarray*}
K_2(\xi_1,\tau_1) &\lesssim& \dfrac{|\xi_1|^{\gamma(s)}}{\langle\sigma_1\rangle^{2b}} \int_{|\eta| \lesssim \langle\sigma_1\rangle+|\xi_1|} \dfrac{d\eta} {|\xi_1|\langle\eta\rangle^{2a}}\\
&\lesssim&\dfrac{|\xi_1|^{\gamma(s)-2a}}{\langle\sigma_1\rangle^{2b}} + \dfrac{|\xi_1|^{\gamma(s)-1}}{\langle\sigma_1\rangle^{2b+2a-1}}
\lesssim 1
\end{eqnarray*}
for $s>-1/4$, $b>1/2$ and $0<a<1/2$ such that $\gamma(s)\leq \min\{1,2a\}$=2a.\\

In the region $B_{4,2}$, by the algebraic relation (\ref{AR2}) we have
\begin{equation*}
\langle\xi_1\rangle\lesssim|\xi_1|\lesssim|\xi_1\xi|\lesssim \langle\tau_1-\xi_1^2\rangle.
\end{equation*} 

Moreover, the change of variables $\eta=\tau_1-\xi_1^2+2\xi\xi_1$, and the restriction in the region $B_{4,2}$ and (\ref{AR2}) give
\begin{equation*}
|\eta|\lesssim\langle\sigma_1\rangle.
\end{equation*} 

Therefore, 
\begin{eqnarray*}
K_2(\xi_1,\tau_1) &\lesssim& \dfrac{\langle\xi_1\rangle^{\gamma(s)}}{\langle\sigma_1\rangle^{2b}} \int_{|\eta| \lesssim \langle\sigma_1\rangle} \dfrac{d\eta} {|\xi_1|\langle\eta\rangle^{2a}}\\
&\lesssim&\dfrac{|\xi_1|^{\gamma(s)-1}}{\langle\sigma_1\rangle^{2b+2a-1}}
\lesssim 1
\end{eqnarray*}
for $s>-1/4$, $b>1/2$ and $0<a<1/2$ such that $\gamma(s)\leq 1$.\\

Finally, we estimate $K_3(\xi_1,\tau_1)$. In the region $B_{4,3}$ we have by the algebraic relation (\ref{AR2}) that 
\begin{equation*}
\langle\xi_1\rangle\lesssim|\xi_1|\lesssim|\xi_1(\xi_1+\xi_2)|\lesssim \langle\sigma_2\rangle.
\end{equation*} 

Therefore Lemma \ref{l3.1}-(\ref{CI2}) implies that
\begin{eqnarray*}
K_3(\xi_1,\tau_1) &\lesssim& \langle\xi_1\rangle^{\gamma(s)-2b}\int \dfrac{1}{\langle\tau_2+\xi_2^2+2\xi_1^2+2\xi_1\xi_2\rangle^{2a}}d\xi_1 \\
&\lesssim&  1
\end{eqnarray*}
for $a>1/4$, $b>1/2$ and $s>-1/4$ which implies $\gamma(s)\leq 2b$.\\
\fim

We finish this section with a result that will be useful in the proof of Theorem \ref{t1.3}. 
\begin{coro}\label{c3.3}
Let $s >-1/4$ and $a, b\in\R$ given in Theorem \ref{t1.1}. For $s'> s$ we have
\begin{equation}\label{BE2}
\left\|
\left(\dfrac{|\xi|^2\widetilde{uv}(\tau,\xi)}{2i\gamma(\xi)}\right)^{\sim^{-1}}
\right\|_{X_{s',-a}}\leq c\left\|u\right\|_{X_{s',b}}\left\|v\right\|_{X_{s,b}} + c\left\|u\right\|_{X_{s,b}}\left\|v\right\|_{X_{s',b}}.
\end{equation}
\end{coro}

\textbf{Proof. } The result is a direct consequence of Theorem \ref{t1.1} and the inequality
\begin{equation*}
 \langle\xi\rangle^{s'}\leq  \langle\xi\rangle^{s}\langle\xi_1\rangle^{s'-s}  +\langle\xi\rangle^{s}\langle\xi-\xi_1\rangle^{s'-s} .
\end{equation*}
\fim
\section{Counterexample to the bilinear estimates (\ref{BE})}

\textbf{Proof of Theorem \ref{t1.2}.} Let $A_N$ denote the set
\begin{equation*}
A_N=\left\{  
\begin{array}{l l}
(\tau,\xi)\in \R^2:(\tau,\xi)=(N^2,N)+\alpha\vec{\eta}+\beta\vec{\gamma}\\
0\leq\alpha\leq N,\peq 0\leq\beta\leq N^{-1},\\
\vec{\eta}=\dfrac{(2N,1)}{\sqrt{1+4N^2}},\peq \vec{\gamma}=\dfrac{(-1,2N)}{\sqrt{1+4N^2}}
\end{array} \right\}
\end{equation*}

and define $f_N(\tau,\xi)=\chi_{A_{N}}$, $g_N(\tau,\xi)=\chi_{-A_{N}}$ where 
\begin{equation*}
-A_N=\left\{  
(\tau,\xi)\in \R^2:-(\tau,\xi)\in A_N\right\}.
\end{equation*}

It is easy to see that 
\begin{equation}\label{N1}
\|f_N\|_{L^2_{\tau,\xi}}=\|g_N\|_{L^2_{\tau,\xi}}=1.
\end{equation} 

Now, let $u_N, v_N\in X_{s,b}$ such that $f_N(\tau,\xi)\equiv \langle|\tau|-\xi^2\rangle^b\langle\xi\rangle^{s} \widetilde{u}_N(\tau,\xi)$ and $g_N(\tau,\xi)\equiv \langle|\tau|-\xi^2\rangle^b\langle\xi\rangle^{s} \widetilde{v}_N(\tau,\xi)$.\\

Therefore, from Lemma \ref{l3.3}-(\ref{LN}) and the fact that 
\begin{eqnarray*}
||\tau|-\xi^2|\leq \min\{|\tau-\xi^2|,|\tau+\xi^2|\}
\end{eqnarray*}
we obtain
\begin{small}
\begin{eqnarray*}
\left\|
\left(\dfrac{|\xi|^2\widetilde{u_Nv_N}(\tau,\xi)}{2i\gamma(\xi)}\right)^{\sim^{-1}}
\right\|_{X_{s,-a}}\equiv
\end{eqnarray*}
\begin{eqnarray*}
&\equiv& \left\|\dfrac{|\xi|^2\langle\xi\rangle^{s}} {\gamma(\xi)\langle|\tau|-\xi^2\rangle^{a}}
\int\!\!\!\!\!\int 
\dfrac{ f_N(\tau_1,\xi_1)\langle\xi_1\rangle^{-s} g_N(\tau-\tau_1,\xi-\xi_1)\langle\xi-\xi_1\rangle^{-s} d\tau_1 d\xi_1}{\langle|\tau-\tau_1|-\gamma(\xi-\xi_1)\rangle^b \langle|\tau_1|-\gamma(\xi_1)\rangle^b}\right\|_{L^2_{\tau,\xi}}\\
&\gtrsim&\left\|\dfrac{|\xi|^2\langle\xi\rangle^{s}} {\gamma(\xi)\langle\tau-\xi^2\rangle^{a}}
\int\!\!\!\!\!\int 
\dfrac{ f_N(\tau_1,\xi_1)\langle\xi_1\rangle^{-s} g_N(\tau-\tau_1,\xi-\xi_1)\langle\xi-\xi_1\rangle^{-s} d\tau_1 d\xi_1}{\langle\tau-\tau_1+(\xi-\xi_1)^2\rangle^b \langle\tau_1-\xi_1^2\rangle^b}\right\|_{L^2_{\tau,\xi}}\\
&\equiv& B_N
\end{eqnarray*}
\end{small}

From the definition of $A_N$ we have
\begin{itemize}
\item [(i)] If $(\tau_1,\xi_1)\in \textrm{supp}\,f_N$ and $(\tau-\tau_1,\xi-\xi_1)\in \textrm{supp}\,g_N$ then there exists $c>0$ such that
\begin{equation*}
|\tau_1-\xi_1^2|\leq c \ppeq \textrm{ and } \ppeq |\tau-\tau_1+(\xi-\xi_1)^2|\leq c. 
\end{equation*}
\item [(ii)]$f\ast g(\tau,\xi)\geq \chi_{R_N}(\tau,\xi)$,

where $R_N$ is the rectangle of dimensions $cN\times\left({cN}\right)^{-1}$ with one of the vertices at the origin and the longest side pointing in the $(1,2N)$ direction.
\item [(iii)] There exists a positive constant $c>0$ such that 
\begin{equation*}
 N\leq \xi_1\leq N+c,\ppeq N\leq \xi_1-\xi\leq N+c
\end{equation*}
and, therefore $|\xi|\sim c$.
\end{itemize}
Moreover, combining the following algebraic relation
\begin{eqnarray*}
(\tau-\tau_1+(\xi-\xi_1)^2)+(\tau_1-\xi_1^2)-(\tau-\xi^2)=2\xi(\xi_1-\xi)
\end{eqnarray*}
with (i) and (iii) we obtain
\begin{equation}\label{TAUXI}
|\tau-\xi^2|\lesssim N.
\end{equation}

Therefore $(\ref{N1})$, (i), (ii), (iii) and $(\ref{TAUXI})$ imply that
\begin{eqnarray*}
1\peq \gtrsim\peq B_N&\gtrsim& \dfrac{N^{-2s}}{N^{a}}\left\|\dfrac{|\xi|^2}{\gamma(\xi)}\chi_{R_N}\right\| _{L^2_{\tau,\xi}}\\
&\gtrsim&\dfrac{N^{-2s}}{N^{a}}\left(\int\!\!\!\!\int_{\{|\xi|\geq 1/2\}} \chi^2_{R_N}(\tau,\xi) \right)^{1/2} \gtrsim N^{-2s-a}.
\end{eqnarray*}

Letting $N\rightarrow\infty$, this inequality is possible only when $-2s-a\leq 0$ which yields the result since $a<1/2$.\\
\fim


\section{Local Well-posedness}

\textbf{Proof of Theorem \ref{t1.3}.} 
\begin{enumerate}
\item \textit{Existence.}\\

 For $({\phi},{\psi})\in H^s(\R)\times H^{s-1}(\R)$, with $s>-1/4$, and $T\leq 1$ we define the integral equation
\begin{equation}\label{IE}
\Gamma_T(u)(t)= \theta(t)\left(V_c(t)\phi+V_s(t)\psi_x\right)+\theta_T(t)\int_{0}^{t} V_s(t-t')(u^2)_{xx}(t')dt'.
\end{equation}

Our goal is to use the Picard fixed point theorem to find a solution
\begin{equation*}
\Gamma_T(u)=u.
\end{equation*}

Let $s>-1/4$ and $a,b\in \R$ such that Theorem \ref{t1.1} holds, that is, $1/4<a<1/2<b$ and $1-(a+b)\equiv \delta >0$.\\

Therefore using (\ref{LP}), Lemma \ref{L22}-($ii$) with $b'=-a$ and (\ref{BE}) we obtain
\begin{eqnarray}\label{C1}
\begin{split}
\|\Gamma_T(u)\|_{X_{s,b}}
&\leq c\left(\|\phi\|_{H^s}+ \|\psi\|_{H^{s-1}}+T^{\delta}\left\|u\right\|_{X_{s,b}}^2\right)\\
\|\Gamma_T(u)-\Gamma_T(v)\|_{X_{s,b}}&\leq cT^{\delta} \left\|u+v\right\|_{X_{s,b}}\left\|u-v\right\|_{X_{s,b}}.
\end{split}
\end{eqnarray}

We define 
\begin{equation*}
X_{s,b}(d)=\left\{u\in X_{s,b}:\|u\|_{X_{s,b}}\leq d\right\}
\end{equation*}
where $d=2c\left(\|\phi\|_{H^s}+\|\psi\|_{H^{s-1}}\right)$.\\

Then choosing 
\begin{equation*}
0<T<\min\left\{\dfrac{1}{(4ca)^{1/\delta}},1\right\}
\end{equation*}
we have that $\Gamma_T:X_{s,b}(d)\rightarrow X_{s,b}(d)$ is a contraction and therefore there exists a unique solution ${u}\in X_{s,b}(d)$ of (\ref{IE}).

Moreover, by Lemma \ref{l11}, we have that $\tilde{u}=u|_{[0,T]} \in C([0,T]:H^s)\cap X_{s,b}^T$ is a solution of (\ref{INT}) in $[0,T]$.
 
\item \textit{If $s'>s$, the result holds in the time interval $[0,T]$ with\\ $T=T(\|\phi\|_{H^s},\|\psi\|_{H^{s-1}})$.}\\

Let $s>-1/4$ and $a,b\in \R$ given in Theorem \ref{t1.1}. For $s'>s$ we consider the closed ball in the Banach space 
\begin{equation*}
W=\left\{u\in X_{s',b}:\|u\|_{s'}=\|u\|_{X_{s,b}}+ \beta\|u\|_{X_{s',b}}< +\infty\right\}
\end{equation*}
where $\beta=\dfrac{\|\phi\|_{H^s}+\|\psi\|_{H^{s-1}}} {\|\phi\|_{H^{s'}}+\|\psi\|_{H^{s'-1}}}$.\\


In view of estimate (\ref{C1}) we obtain
\begin{eqnarray*}
\|\Gamma_T(u)\|_{X_{s,b}}\leq c\left(\|\phi\|_{H^s}+ \|\psi\|_{H^{s-1}}+T^{\delta}\left\|u\right\|_{X_{s,b}}^2\right).
\end{eqnarray*}

Now by Corollary \ref{c3.3} we have
\begin{eqnarray*}
\|\Gamma_T(u)\|_{X_{s',b}}
&\leq& c\left(\|\phi\|_{H^{s'}}+ \|\psi\|_{H^{s'-1}}+T^{\delta}\left\|u\right\|_{X_{s',b}} \left\|u\right\|_{X_{s,b}}\right)\\
&\leq& \dfrac{c}{\beta}\left(\|\phi\|_{H^{s}}+ \|\psi\|_{H^{s-1}}+T^{\delta}\left\|u\right\|_{s'}^2\right).
\end{eqnarray*}

Therefore
\begin{eqnarray*}
\|\Gamma_T(u)\|_{s'}&\leq 2c\left(\|\phi\|_{H^{s}}+ \|\psi\|_{H^{s-1}}+T^{\delta}\left\|u\right\|_{s'}^2\right).
\end{eqnarray*}

The same argument gives
\begin{eqnarray*}
\|\Gamma_T(u)-\Gamma_T(v)\|_{s'}\leq 2cT^{\delta}\left\|u+v\right\|_{s'} \left\|u-v\right\|_{s'}.
\end{eqnarray*}

Then we define in $W$ the closed ball centered at the origin with radius $d'=4c\left(\|\phi\|_{H^{s}}+ \|\psi\|_{H^{s-1}}\right)$ and choose
\begin{equation*}
0<T<\min\left\{\dfrac{1}{(8cd')^{1/\delta}},1\right\}.
\end{equation*}

Thus we have that $F_T$ is a contraction and therefore there exists a solution with $T=T(\|\phi\|_{H^s},\|\psi\|_{H^{s-1}})$.

\item \textit{Uniqueness.} By the fixed point argument used in item 1 we have uniqueness of the solution of the truncated integral equation (\ref{IE}) in the set $X_{s,b}(d)$. We use an argument due to Bekiranov, Ogawa and Ponce \cite{BOP} to obtain the uniqueness of the integral equation (\ref{INT}) in the whole space $X_{s,b}^T$.

Let $T>0$, $u\in X_{s,b}$ be the solution of the truncated integral equation (\ref{IE}) obtained above and $\widetilde{v}\in X_{s,b}^T$ be a solution of the integral equation (\ref{INT}) with the same initial data. Fix an extension $v \in X_{s,b}$, therefore, for some $T^{\ast}<T<1$ to be fixed later, we have
\begin{equation*}
v(t)= \theta(t)\left(V_c(t)\phi+V_s(t)\psi_x\right)+\theta_T(t)\int_{0}^{t} V_s(t-t')(v^2)_{xx}(t')dt'
\end{equation*}
for all $t\in[0,T^{\ast}]$.

Fix $M\geq0$ such that
\begin{equation}\label{COTAM}
\max\left\{\|u\|_{X_{s,b}},\|v\|_{X_{s,b}}\right\}\leq M.
\end{equation}

Taking the difference $u-v$, by definition of $X_{s,b}^{T^{\ast}}$, we have that for any $\e>0$, there exists $w\in X_{s,b}$ such that for all $t\in[0,T^{\ast}]$
\begin{equation*}
w(t)=u(t)-v(t)
\end{equation*}
and
\begin{equation}\label{DEFX}
\|w\|_{X_{s,b}} \leq \|u-v\|_{X_{s,b}^{T^{\ast}}}+\varepsilon.
\end{equation}

Define
\begin{equation*}
\widetilde{w}(t)=\theta_T(t)\int_{0}^{t} V_s(t-t')(w(t')u(t')+w(t')v(t'))_{xx}(t')dt'.
\end{equation*}

We have that, $\widetilde{w}(t)=u(t)-v(t)$, for all $t\in[0,T^{\ast}]$. Therefore, from Definition \ref{BL}, Lemma \ref{L22}-($ii$), (\ref{BE}) and (\ref{COTAM}) it follows that
\begin{eqnarray}\label{UNI}
 \|u-v\|_{X_{s,b}^{T^{\ast}}}\leq \|\widetilde{w}\|_{X_{s,b}}
\leq 2cMT^{\ast\delta} \left\|w\right\|_{X_{s,b}}.
\end{eqnarray}

Choosing $T^{\ast}>0$ such that $2cMT^{\ast\delta}\leq 1/2$, by (\ref{DEFX}) and (\ref{UNI}), we have
\begin{equation*}
\|u-v\|_{X_{s,b}^{T^{\ast}}}\leq 2\e.
\end{equation*}

Therefore $u=v$ on $[0,T^{\ast}]$. Now, since the argument does not depend on the initial data, we can iterate this process a finite number of times to extend the uniqueness result in the whole existence interval $[0,T]$.

\item \textit{Map data-solution is locally Lipschitz.} Combining an identical argument to the one used in the existence proof with Lemma \ref{l11}, one can easily show that the map data-solution is locally Lipschitz.
\end{enumerate}
\fim

\section{Proof of Theorems \ref{t4.1}-\ref{t4.2}}

\textbf{Proof of Theorem \ref{t4.1}} Suppose that there exists a space $X_T$ satisfying the conditions of the theorem for $s<-2$ and $T>0$. Let $\phi,\rho \in H^{s}(\R)$ and define $u(t)=V_c(t)\phi$, $v(t)=V_c(t)\rho$. In view of (\ref{i}), (\ref{ii}), (\ref{iii}) it is easy to see that the following inequality must hold
\begin{equation}\label{iv}
\sup_{1\leq t\leq T}\left\|\int_{0}^{t}V_s(t-t')(V_c(t')\phi V_c(t')\rho)_{xx}(t')dt'\right\|_{H^{s}(\R)}\leq c\left\|\phi\right\|_{H^{s}(\R)}\left\|\rho\right\|_{H^{s}(\R)}.
\end{equation}

We will see that (\ref{iv}) fails for an appropriate choice of $\phi$, $\rho$, which would lead to a contradiction.

Define
\begin{equation*}
\widehat{\phi}(\xi)=N^{-s}\chi_{[-N,-N+1]} \textrm{\para and \para} \widehat{\rho}(\xi)=N^{-s}\chi_{[N+1,N+2]},
\end{equation*} 
where $\chi_A(\cdot)$ denotes the characteristic function of the set $A$.\\

We have
\begin{equation*}
\left\|\phi\right\|_{H^s(\R)}, \left\|\rho\right\|_{H^s(\R)} \sim 1.
\end{equation*} 

Recall that $\gamma(\xi)\equiv \sqrt{{\xi}^2+{\xi}^4}$. By the definitions of $V_c$, $V_s$ and Fubini's Theorem, we have
\begin{equation*}
\left(\int_{0}^{t}V_s(t-t')(V_c(t')\phi V_c(t')\rho)_{xx}(t')dt'\right)^{\wedge_{(x)}}(\xi)=
\end{equation*}
\begin{eqnarray*}
&=&\int_{-\infty}^{+\infty}-\dfrac{|\xi|^2}{8i\gamma(\xi)} \widehat{\phi}(\xi-\xi_1)\widehat{\rho}(\xi_1)K(t,\xi,\xi_1)d\xi_1\\
&=&\int_{A_{\xi}}-\dfrac{|\xi|^2}{8i\gamma(\xi)} N^{-2s}K(t,\xi,\xi_1) d\xi_1
\end{eqnarray*}
where
\begin{equation*}
A_{\xi}=\left\{\xi_1: \xi_1\in \textrm{supp}(\widehat{\rho}) \textrm{ and } \xi-\xi_1\in \textrm{supp}(\widehat{\phi})\right\}
\end{equation*}
and
\begin{equation*}
K(t,\xi,\xi_1)\equiv \int_{0}^{t}\sin((t-t')\gamma(\xi))\cos(t'\gamma(\xi-\xi_1)) \cos(t'\gamma(\xi_1)) dt'.
\end{equation*}

Note that for all $\xi_1\in \textrm{supp}(\widehat{\rho})$ and $\xi-\xi_1\in \textrm{supp}(\widehat{\phi})$ we have
\begin{equation*}
\gamma(\xi-\xi_1),\gamma(\xi_1)\sim N^2 \textrm{ and } 1\leq \xi \leq 3.
\end{equation*}

On the other hand, since $s<-2$, we can choose $\varepsilon>0$ such that
\begin{equation}\label{CON}
-2s-4-2\varepsilon>0.
\end{equation}

Let $t=\dfrac{1}{N^{2+\varepsilon}}$, then for $N$ sufficiently large we have
\begin{equation*}
\cos(t'\gamma(\xi-\xi_1)), \cos(t'\gamma(\xi_1))\geq 1/2
\end{equation*}
and
\begin{equation*}
\sin((t-t')\gamma(\xi))\geq c(t-t')\gamma(\xi),
\end{equation*}
for all $0\leq t'\leq t$, $1\leq \xi \leq 3$ and $\xi_1\in \textrm{supp}(\widehat{\eta})$.\\

Therefore
\begin{equation*}
K(t,\xi,\xi_1)\gtrsim \int_{0}^{t}(t-t')\gamma(\xi)dt'\gtrsim \gamma(\xi)\dfrac{1}{2N^{4+2\varepsilon}}.
\end{equation*}

For $3/2\leq\xi\leq 5/2$ we have that $\textrm{mes}(A_{\xi})\gtrsim 1$. Thus, from (\ref{iv}) we obtain 
\begin{eqnarray*}
1 &\gtrsim& \sup_{1\leq t\leq T}\left\|\int_{0}^{t}V_s(t-t')(V_c(t')\phi V_c(t')\rho)_{xx}(t')dt'\right\|_{H^{s}(\R)}\\
&\gtrsim& \sup_{1\leq t\leq T}\left(\int_{3/2}^{5/2}\left(1+|\xi|^2\right)^s\left| \int_{A_{\xi}}\dfrac{|\xi|^2}{8i\gamma(\xi)}N^{-2s} K(t,\xi,\xi_1) d\xi_1\right|^2d\xi\right)^{1/2}\\
&\gtrsim& N^{-2s-4-2\varepsilon}, \textrm{ for all } N\gg 1
\end{eqnarray*}
which is in contradiction with (\ref{CON}).\\
\fim

\textbf{Proof of Theorem \ref{t4.2}} Let $s<-2$ and suppose that there exists $T>0$ such that the flow-map $S$ defined in (\ref{DM}) is $C^2$. When $(\phi,\psi)\in H^{s}(\R)\times H^{s-1}(\R)$, we denote by $u_{(\phi,\psi)}\equiv S(\phi,\psi)$ the solution of the IVP (\ref{NLB}) with $f(u)=u^2$, $u_0=\phi$ and $u_1=\psi_x$, that is 
\begin{equation*}
u_{(\phi,\psi)}(t)= V_c(t)\phi+V_s(t)\psi_x+\int_{0}^{t}V_s(t-t')(u_{(\phi,\psi)}^2)_{xx}(t')dt'.
\end{equation*}

The Fr\'echet derivative of $S$ at $(\omega,\zeta)$ in the direction $(\phi,\bar{\phi})$ is given by
\begin{eqnarray}\label{FRECHET}
d_{(\phi,\bar{\phi})}S(\omega,\zeta)= V_c(t)\phi+V_s(t)\bar{\phi}_x+2\int_{0}^{t}V_s(t-t')(u_{(\phi,\psi)}(t') d_{(\phi,\bar{\phi})}S(\omega,\zeta)(t'))_{xx}dt'.
\end{eqnarray}

Using the well-posedness assumption we know that the only solution for initial data $(0,0)$ is $u_{(0,0)}\equiv S(0,0)=0$. Therefore, (\ref{FRECHET}) yields
\begin{eqnarray*}
d_{(\phi,\bar{\phi})}S(0,0)=V_c(t)\phi+V_s(t)\bar{\phi}_x.
\end{eqnarray*}

Computing the second Fr\'echet derivative at the origin in the direction $((\phi,\bar{\phi}),(\rho,\bar{\rho}))$, we obtain
\begin{eqnarray*}
d^2_{(\phi,\bar{\phi}),(\rho,\bar{\rho})}S(0,0)=
\end{eqnarray*}
\begin{eqnarray*}
=2\int_{0}^{t}V_s(t-t')\left[(V_c(t')\phi+V_s(t')\bar{\phi}_x) (V_c(t')\rho+V_s(t')\bar{\rho}_x)\right]_{xx}dt'.
\end{eqnarray*}

Taking $\bar{\phi},\bar{\rho}=0$, the assumption of $C^2$ regularity of $S$ yields
\begin{equation*}
\sup_{1\leq t\leq T}\left\|\int_{0}^{t}V_s(t-t')(V_c(t')\phi V_c(t')\rho)_{xx}(t')dt'\right\|_{H^{s}(\R)}\leq c\left\|\phi\right\|_{H^{s}(\R)}\left\|\rho\right\|_{H^{s}(\R)}
\end{equation*}
which has been shown to fail in the proof of Theorem \ref{t4.1}.\\
\fim

\centerline{\textbf{Acknowledgment}}

This paper is prepared under the guidance of my advisor Felipe Linares (IMPA). I want to take the opportunity to express my sincere gratitude to him. I also thank Aniura Milan\'es (UFMG), Dami\'an Fern\'andez (Unicamp) and Didier Pilod (UFRJ) for fruitful conversations concerning this work. The author was partially supported by CNPq-Brazil.


\begin{thebibliography}{10}

\bibitem{BT}
I.~Bejenaru and T.~Tao.
\newblock Sharp well-posedness and ill-posedness results for a quadratic
  non-linear {S}chr\"odinger equation.
\newblock {\em J. Funct. Anal.}, 233(1):228--259, 2006.

\bibitem{BOP}
D.~Bekiranov, T.~Ogawa, and G.~Ponce.
\newblock Interaction equations for short and long dispersive waves.
\newblock {\em J. Funct. Anal.}, 158(2):357--388, 1998.

\bibitem{BS}
J.~L. Bona and R.~L. Sachs.
\newblock Global existence of smooth solutions and stability of solitary waves
  for a generalized {B}oussinesq equation.
\newblock {\em Comm. Math. Phys.}, 118(1):15--29, 1988.

\bibitem{B}
J.~Bourgain.
\newblock Fourier transform restriction phenomena for certain lattice subsets
  and applications to nonlinear evolution equations. {I} and {II}. {T}he
  {K}d{V}-equation.
\newblock {\em Geom. Funct. Anal.}, 3(3):107--156, 209--262, 1993.

\bibitem{B2}
J.~Bourgain.
\newblock Periodic {K}orteweg de {V}ries equation with measures as initial
  data.
\newblock {\em Selecta Math. (N.S.)}, 3(2):115--159, 1997.

\bibitem{BOU}
J.~Boussinesq.
\newblock Th\'eorie des ondes et des remous qui se propagent le long d'un canal
  rectangulaire horizontal, en communiquant au liquide continu dans 21 ce canal
  des vitesses sensiblement pareilles de la surface au fond.
\newblock {\em J. Math. Pures Appl.}, 17(2):55--108, 1872.

\bibitem{FLS}
F.~Falk, E.~Laedke, and K.~Spatschek.
\newblock Stability of solitary-wave pulses in shape-memory alloys.
\newblock {\em Phys. Rev. B}, 36(6):3031--3041, 1987.

\bibitem{FG}
Y.-F. Fang and M.~G. Grillakis.
\newblock Existence and uniqueness for {B}oussinesq type equations on a circle.
\newblock {\em Comm. Partial Differential Equations}, 21(7-8):1253--1277, 1996.

\bibitem{G}
J.~Ginibre.
\newblock Le probl\`eme de {C}auchy pour des {EDP} semi-lin\'eaires
  p\'eriodiques en variables d'espace (d'apr\`es {B}ourgain).
\newblock {\em Ast\'erisque}, (237):Exp.\ No.\ 796, 4, 163--187, 1996.
\newblock S\'eminaire Bourbaki, Vol.\ 1994/95.

\bibitem{GTV}
J.~Ginibre, Y.~Tsutsumi, and G.~Velo.
\newblock On the {C}auchy problem for the {Z}akharov system.
\newblock {\em J. Funct. Anal.}, 151(2):384--436, 1997.

\bibitem{KPV1}
C.~E. Kenig, G.~Ponce, and L.~Vega.
\newblock A bilinear estimate with applications to the {K}d{V} equation.
\newblock {\em J. Amer. Math. Soc.}, 9(2):573--603, 1996.

\bibitem{KPV2}
C.~E. Kenig, G.~Ponce, and L.~Vega.
\newblock Quadratic forms for the {$1$}-{D} semilinear {S}chr\"odinger
  equation.
\newblock {\em Trans. Amer. Math. Soc.}, 348(8):3323--3353, 1996.

\bibitem{FL}
F.~Linares.
\newblock Global existence of small solutions for a generalized {B}oussinesq
  equation.
\newblock {\em J. Differential Equations}, 106(2):257--293, 1993.

\bibitem{MST1}
L.~Molinet, J.~C. Saut, and N.~Tzvetkov.
\newblock Ill-posedness issues for the {B}enjamin-{O}no and related equations.
\newblock {\em SIAM J. Math. Anal.}, 33(4):982--988 (electronic), 2001.

\bibitem{MST2}
L.~Molinet, J.-C. Saut, and N.~Tzvetkov.
\newblock Well-posedness and ill-posedness results for the
  {K}adomtsev-{P}etviashvili-{I} equation.
\newblock {\em Duke Math. J.}, 115(2):353--384, 2002.

\bibitem{TM}
M.~Tsutsumi and T.~Matahashi.
\newblock On the {C}auchy problem for the {B}oussinesq type equation.
\newblock {\em Math. Japon.}, 36(2):371--379, 1991.

\bibitem{T}
N.~Tzvetkov.
\newblock Remark on the local ill-posedness for {K}d{V} equation.
\newblock {\em C. R. Acad. Sci. Paris S\'er. I Math.}, 329(12):1043--1047,
  1999.

\bibitem{Z}
V.~Zakharov.
\newblock On stochastization of one-dimensional chains of nonlinear
  oscillators.
\newblock {\em Sov. Phys. JETP}, 38:108--110, 1974.

\end{thebibliography}

E-mail: farah@impa.br
\end{document}